\documentclass[12pt]{article}
\usepackage{amssymb,latexsym,pstcol,pst-plot}

\newcommand{\ben}{\begin{enumerate}}
\newcommand{\een}{\end{enumerate}}
\newcommand{\ble}{\begin{lem}}
\newcommand{\ele}{\end{lem}}
\newcommand{\bth}{\begin{thm}}
\renewcommand{\eth}{\end{thm}}
\newcommand{\bpr}{\begin{prop}}
\newcommand{\epr}{\end{prop}}
\newcommand{\bco}{\begin{cor}}
\newcommand{\eco}{\end{cor}}
\newcommand{\bcon}{\begin{conj}}
\newcommand{\econ}{\end{conj}}
\newcommand{\bde}{\begin{defn}}
\newcommand{\ede}{\end{defn}}
\newcommand{\bex}{\begin{exa}}
\newcommand{\eex}{\end{exa}}
\newcommand{\barr}{\begin{array}}
\newcommand{\earr}{\end{array}}
\newcommand{\btab}{\begin{tabular}}
\newcommand{\etab}{\end{tabular}}
\newcommand{\beq}{\begin{equation}}
\newcommand{\eeq}{\end{equation}}
\newcommand{\bea}{\begin{eqnarray*}}
\newcommand{\eea}{\end{eqnarray*}}
\newcommand{\bce}{\begin{center}}
\newcommand{\ece}{\end{center}}
\newcommand{\bpi}{\begin{picture}}
\newcommand{\epi}{\end{picture}}
\newcommand{\bfi}{\begin{figure} \begin{center}}
\newcommand{\efi}{\end{center} \end{figure}}
\newcommand{\capt}{\caption}
\newcommand{\bsl}{\begin{slide}{}}
\newcommand{\esl}{\end{slide}}

\newcommand{\bib}{thebibliography}
\newcommand{\pf}{{\bf Proof}\hspace{7pt}}

\newcommand{\qqed}{\qquad\rule{1ex}{1ex}}
\newcommand{\Qqed}{\qquad\rule{1ex}{1ex}\medskip}

\newcommand{\hs}[1]{\hspace{#1}}
\newcommand{\hso}[1]{\hspace{-1pt}}

\newcommand{\emp}{\emptyset}

\newcommand{\sbe}{\subseteq}

\newcommand{\iso}{\cong}
\newcommand{\Cong}{\equiv}

\newcommand{\case}[4]{\left\{\barr{ll}#1&\mbox{#2}\\#3&\mbox{#4}\earr\right.}

\def\<{\langle}
\def\>{\rangle}

\newcommand{\ree}[1]{(\ref{#1})}

\newcommand{\ga}{\gamma}
\newcommand{\de}{\delta}

\newcommand{\om}{\omega}

\newcommand{\De}{\Delta}

\newcommand{\ba}{{\bf a}}
\newcommand{\bb}{{\bf b}}
\newcommand{\bc}{{\bf c}}

\newcommand{\bs}{{\bf s}}
\newcommand{\bt}{{\bf t}}

\newcommand{\bbN}{{\mathbb N}}

\newcommand{\cO}{{\cal O}}

\newcommand{\cT}{{\cal T}}

\newcommand{\fS}{{\mathfrak S}}

\newcommand{\Mod}{\mathop{\rm mod}\nolimits}

\newcommand{\dil}{\displaystyle}

\newcommand{\Gca}{\put(20,0){\circle*{3}}}

\newcommand{\Gga}{\put(60,0){\circle*{3}}}

\newcommand{\Gac}{\put(0,20){\circle*{3}}}

\newcommand{\Gcc}{\put(20,20){\circle*{3}}}

\newcommand{\Gec}{\put(40,20){\circle*{3}}}

\newcommand{\Gae}{\put(0,40){\circle*{3}}}

\newcommand{\Gce}{\put(20,40){\circle*{3}}}

\newcommand{\Gee}{\put(40,40){\circle*{3}}}

\newcommand{\Gge}{\put(60,40){\circle*{3}}}

\newcommand{\Gcg}{\put(20,60){\circle*{3}}}

\newcommand{\Geg}{\put(40,60){\circle*{3}}}

\newcommand{\Gcacc}{\put(20,0){\line(0,1){20}}}

\newcommand{\Gcaec}{\put(20,0){\line(1,1){20}}}

\newcommand{\Ggaec}{\put(60,0){\line(-1,1){20}}}

\newcommand{\Gacae}{\put(0,20){\line(0,1){20}}}

\newcommand{\Gacce}{\put(0,20){\line(1,1){20}}}

\newcommand{\Gccce}{\put(20,20){\line(0,1){20}}}

\newcommand{\Gecce}{\put(40,20){\line(-1,1){20}}}

\newcommand{\Gecee}{\put(40,20){\line(0,1){20}}}

\newcommand{\Gaecg}{\put(0,40){\line(1,1){20}}}

\newcommand{\Gcecg}{\put(20,40){\line(0,1){20}}}

\newcommand{\Gceeg}{\put(20,40){\line(1,1){20}}}

\newcommand{\Geecg}{\put(40,40){\line(-1,1){20}}}

\newcommand{\Ggeeg}{\put(60,40){\line(-1,1){20}}}

\newcommand{\Com}{\mathop{\rm Com}\nolimits}

\newtheorem{thm}{Theorem}[section]
\newtheorem{prop}[thm]{Proposition}
\newtheorem{cor}[thm]{Corollary}
\newtheorem{lem}[thm]{Lemma}
\newtheorem{conj}[thm]{Conjecture}
\newtheorem{exa}[thm]{Example}

\begin{document}
\pagestyle{plain}

\title{Congruences for Catalan and Motzkin numbers and related sequences
}
\author{
Emeric Deutsch\\[-5pt]
\small Department of Mathematics\\[-5pt]
\small Polytechnic University\\[-5pt]
\small Brooklyn, NY 11201, USA\\[-5pt]
\small \texttt{deutsch@duke.poly.edu}\\[5pt]
and\\[5pt]
Bruce E. Sagan\\[-5pt]
\small Department of Mathematics\\[-5pt] 
\small Michigan State University\\[-5pt]
\small East Lansing, MI 48824-1027, USA\\[-5pt]
\small \texttt{sagan@math.msu.edu}
}

\date{\today}
\maketitle

\begin{abstract}
We prove various congruences for Catalan and Motzkin numbers as well
as related sequences.  The common thread is that all these sequences can
be expressed in terms of binomial coefficients.  Our techniques are
combinatorial and algebraic: group actions, induction, and
Lucas' congruence for binomial coefficients come into play.  
A number of our results 
settle conjectures of Benoit Cloitre and Reinhard Zumkeller.  The
Thue-Morse sequence appears in several contexts.
\end{abstract}

\section{Introduction}

Let $\bbN$ denote the nonnegative integers.
The divisibility of the {\it Catalan numbers\/}
$$
C_n=\frac{1}{n+1}{2n\choose n},\qquad n\in\bbN,
$$
by primes and prime powers has been completely determined by
Alter and Kubota~\cite{ak:ppp} using arithmetic techniques.  In
particular, the fact that $C_n$ is odd precisely when $n=2^h-1$ for
some $h\in\bbN$ has attracted the attention of several
authors including Deutsch~\cite{deu:idp},
E\u{g}ecio\u{g}lu~\cite{ege:pcn}, and Simion and Ullman~\cite{su:sln}
who found combinatorial explanations of this result.  In the next
section we will derive the theorem which gives
the largest power of 2 dividing any Catalan number by using group
actions.  In addition to its generality, this technique has the advantage
that when $n=2^h-1$ there is 
exactly one fixed point with all the other orbits having size
divisible by 2.  For other congruences which can be proven using the
action of a group, see Sagan's article~\cite{sag:cva}.

By contrast, almost nothing is known about the residues of the {\it Motzkin
numbers\/}
$$
M_n=\sum_{k\ge0}{n\choose 2k} C_k,\qquad n\in\bbN.
$$
In fact, the only two papers dealing with this matter of which we are
aware are the recent articles of Luca~\cite{luc:pfm} about prime
factors of $M_n$ and of Klazar and Luca~\cite{kl:ipm} about
the periodicity of $M_n$ modulo a positive integer.  In
section~\ref{Mn} we will characterize the parity of the Motzkin
numbers as well as three related sequences.  Surprisingly, the
characterizations involve a sequence which 
encodes the lengths of the blocks in the Thue-Morse sequence.  The
block-length sequence was first studied by Allouche et.\
al.~\cite{aab:srt}.  For 
more information about the Thue-Morse sequence in general, the reader
is referred to the survey article of Allouche and Shallit~\cite{as:ups}.

Section~\ref{cbt} is devoted to congruences for the central binomial
and trinomial coefficients.  We are able to use these results to
describe the Motzkin numbers and their relatives modulo 3.
They also prove various conjectures of Benoit
Cloitre~\cite{clo:pc} and Reinhard Zumkeller~\cite{slo:ole}.  The
Thue-Morse sequence appears again.  Our 
main tool in this section is  Lucas' congruence
for multinomial coefficients~\cite{luc:cne}.

Our final section is a collection of miscellaneous results and
conjectures about sequences related to binomial coefficients.  These
include the Ap\'ery numbers, the central Delannoy and Eulerian
numbers, Gould's sequence, and the sequence enumerating noncrossing
graphs.

\section{Catalan numbers}
\label{Cn}

If $n,m\in\bbN$ with $m\ge2$ then the {\it order of $n$ modulo $m$\/}
is
$$
\om_m(n)=\mbox{largest power of $m$ dividing $n$.}
$$
If the base $m$ expansion of $n$ is
\beq
\label{exp}
n=n_0 +n_1 m+n_2 m^2 +\cdots
\eeq
then let
$$
\De_m(n)=\{i\ :\ n_i=1\}
$$
and
$$
\de_m(n) =|\De_m(n)|
$$
where the absolute value signs denotes cardinality. 
We will also use a pound sign for this purpose. If a subscript $m$ is not
used then we are assuming $m=2$ and in this case $\de(n)$ is also the
sum of the digits in the base 2 expansion of $n$.

We wish to prove the following theorem.
\bth
\label{Cnthm}
For $n\in\bbN$ we have
$$
\om(C_n)=\de(n+1)-1.
$$
\eth
Note as an immediate corollary that $C_n$ is odd if and
only if $n=2^h-1$ for some $h\in\bbN$.  It is easy to prove this
theorem from Kummer's result about the order of a binomial
coefficient~\cite{kum:ear} (or see~\cite[pp.\ 270--271]{dic:htn}).
However, we wish to give a combinatorial proof.

We will use a standard interpretation of $C_n$ using binary trees.
A {\it binary tree\/} $T$ is a tree with a root $r$ where every vertex has a left
child, or a right child, or both, or neither.  
Note that this differs from the convention where a vertex in a binary
tree must have no children or both children.
It will also be
convenient to consider $T=\emp$ as a binary tree.  
With this convention, any nonempty tree can be written as $T=(T',T'')$
where $T'$ and $T''$ are the subtrees generated by the left child and
by the right child of $r$, respectively.   (The subtree {\it generated\/}
by a vertex $v$ of $T$ consists of $v$ and all its descendants.)
Let $\cT_n$ be the
set of all binary trees on $n$ vertices.  Then it is well-known that
$|\cT_n|=C_n$ for all $n\in\bbN$.

The {\it height\/} of a vertex $v$ is the length of the unique path
from the root $r$ to $v$.  A {\it complete binary tree\/} $T_h$
has all $2^i$ possible vertices at height $i$ for $0\le i\le h$ and no
other vertices.  Let $G_h$ be the group of automorphisms of $T_h$ as a
rooted tree.  We will need some facts about $G_h$.
\ble
\label{Gh}
We have the following 
\ben
\item[(1)] If $h=0$ then $G_0=\{e\}$ where $e$ is the identity element,
  and if $h\ge1$ then
$$
G_h = Z_2 \wr G_{h-1}
$$
where $Z_2$ is the cyclic group of order 2 and $\wr$ is 
wreath product.
\item[(2)] $\# G_h= 2^{2^h-1}$.
\item[(3)]  If $G_h$ acts on a set and $\cO$ is an orbit of the action then
$\#\cO$ is a power of 2.
\een
\ele
\pf
The proof of (1) follows by noting that $T_h=(T_{h-1},T_{h-1})$ for
$h\ge1$.  Then (2) is an easy induction on $h$ using (1).  Finally,
(3) is a consequence of (2) and the fact that for any group action the size of
an orbit always divides the order of the group.
\Qqed

Now $G_n$ acts on $\cT_n$ in the obvious way.  It is this action which
will permit us to calculate $\om(C_n)$.  Recall the {\it double
factorial\/}
$$
(2d)!!=(2d-1)(2d-3)\cdots 3\cdot 1.
$$
\ble
\label{cO}
For $n\in\bbN$, let $d=\de(n+1)-1$.  Then given any orbit $\cO$ of
$G_n$ acting on $\cT_n$ we have
$$
\om(\#\cO)\ge d
$$
with equality for exactly $(2d)!!$ orbits.
\ele
\pf
We will induct on $n$ with the result being trivial for $n=0$.  For
$n\ge1$ let $T=(T',T'')\in\cT_n$.  We also let $n'$ and $n''$ be the
number of vertices of $T'$ and $T''$ respectively, as well as setting
$d'=\de(n'+1)-1$ and $d''=\de(n''+1)-1$.  Clearly
$n+1=(n'+1)+(n''+1)$.  It follows that 
\beq
\label{d}
d\le d'+d''+1
\eeq
with equality if and only if we have a disjoint union 
$\De(n+1)=\De(n'+1)\uplus\De(n''+1)$.

Let $\cO(T)$ denote the orbit of $T$.  Then 
\beq
\label{cOeq}
|\cO(T)|=
\case{|\cO(T')|^2}{if $T'\iso T''$,}
{2|\cO(T')||\cO(T'')|}{otherwise.\rule{0pt}{20pt}}
\eeq
Also we have, by induction, $\om(\#\cO(T'))\ge d'$ and
$\om(\#\cO(T''))\ge d''$.

First consider the case when $T'\iso T''$.  Then $n'=n''$ and so
equation~\ree{d} gives $d<2d'+1$.  Now from~\ree{cOeq} we obtain
$$
\om(\#\cO(T))=2\om(\#\cO(T'))\ge 2d'\ge d
$$
as desired for the first half of the lemma.  If we actually have 
$\om(\#\cO(T))= d$ then this forces $2d'=d$.  But since $n'=n''$ we
also have $n+1=2(n'+1)$ and so $d=d'$.  This can only happen if
$d=d'=0$ and consequently $n=2^h-1$ for some $h$.  But by the third
part of the previous lemma, $T_h$ is
the unique tree with $2^h-1$ vertices and $\om(\#\cO(T))=0$.  
Since in this case $(2d)!!=0!!=1$, we have proven the
present lemma when $T'\iso T''$.

Now consider what happens when $T'\not\iso T''$.  Using
equations~\ree{d} and~\ree{cOeq} as before gives
$$
\om(\#\cO(T))=\om(\#\cO(T'))+\om(\#\cO(T'))+1\ge d'+d''+1\ge d
$$
and again the first half of the lemma follows.  When 
$\om(\#\cO(T))= d$ then we must have 
$\om(\#\cO(T'))= d'$, $\om(\#\cO(T''))= d''$, and
$\De(n+1)=\De(n'+1)\uplus\De(n''+1)$.  Using~\ree{cOeq} to count orbits
and induction it follows that we will be done if we can show
\beq
\label{2d!!}
(2d)!!=\frac{1}{2}\sum_{k=1}^d {d+1\choose k}(2k-2)!!(2d-2k)!!
\eeq
for $d\ge1$.  Rewriting this equation in hypergeometric series form we
obtain the equivalent identity
$$
\rule{0pt}{5pt}_2 F_1\left(\barr{c} -d-1,\ -1/2\\ 1/2-d\earr;1\right)=0
$$
which is true by Vandermonde's convolution.
\Qqed

We can now prove Theorem~\ref{Cnthm}.  Since the orbits of a group
action partition the set acted on, we can use Lemma~\ref{Gh} (3) and
Lemma~\ref{cO} to write
$$
C_n=\#\cT_n=(2d)!! 2^d+k 2^{d+1}
$$
for some $k\in\bbN$.  Since $(2d)!!$ is odd we can conclude
$\om(C_n)=d=\de(n+1)-1$ as desired.
\Qqed

The reader may not be happy with the last step in the proof of
Lemma~\ref{cO} since its appeal to the theory of hypergeometric
series is not combinatorial.  So we wish to give a bijective proof
of equation~\ree{2d!!}.  For this, we will interpret the double
factorial in terms of binary total partitions, an object introduced
and enumerated by Schr\"oder~\cite{sch:vcp}.  Given a set $S$ then a
{\it binary total partition of $S$\/} is an unordered rooted tree $B$
satisfying the following restrictions.
\ben
\item Every vertex of $B$ has 0 or 2 children.
\item Every vertex of $B$ is labeled with a subset of $S$ in such a
way that
  \ben 
  \item the root is labeled with $S$ and the leaves with the 1-element
  subsets of $S$,
  \item if a vertex is labeled with $A$ and its children with $A',A''$
  then $A=A'\uplus A''$.
  \een
\een

\thicklines
\setlength{\unitlength}{2pt}
\bfi
\bpi(70,70)(-10,-10)
\put(20,-5){\makebox(0,0){$1$}}
\Gca 
\put(60,-5){\makebox(0,0){$4$}}
\Gga
\put(0,15){\makebox(0,0){$3$}}
\Gac 
\put(45,20){\makebox(0,0){$14$}}
\Gec
\put(10,40){\makebox(0,0){$134$}}
\Gce 
\put(60,35){\makebox(0,0){$2$}}
\Gge
\put(40,65){\makebox(0,0){$1234$}}
\Geg
\Gcaec \Ggaec 
\Gacce \Gecce
\Gceeg \Ggeeg
\epi
\capt{A total binary partition}\label{tbp}
\efi

For example, if $S=\{1,2,3,4\}$ then one possible total binary
partition is displayed in Figure~\ref{tbp}.  Let $b_d$ be the number
of total binary partitions on set $S$ with $|S|=d$.  Then 
$$
b_{d+1}=(2d)!!
$$
For proofs of this fact, including a combinatorial one, see the text of 
Stanley~\cite[Example 5.2.6]{sta:ec2}.  

It is now easy to prove~\ree{2d!!} combinatorially.  The left side
counts total binary partitions $B$ of a set $S$ with $|S|=d+1$. For
the right side, note that each
such $B$ can be formed uniquely by writing $S=S'\uplus S''$, letting
$S'$ and $S''$ label the children of the root, and then forming total
binary partitions on $S'$ and $S''$ to create the rest of $B$.  If
$\#S'=k$ then there are ${d+1\choose k}$ choices for $S'$ (after
which, $S''$ is uniquely determined).  The factors $(2k-2)!!$ and
$(2d-2k)!!$  count the number of ways to put total
binary partitions on $S'$ and $S''$, respectively.  Finally, we must
sum over all possible $k$ and divide by 2 since the tree is unordered.
This completes the combinatorial proof of~\ree{2d!!}.

\section{Motzkin numbers and related sequences}
\label{Mn}

To find the parity of $M_n$ we must first introduce a related
sequence.  Define $\bc=(c_0,c_1,c_2,\ldots)=(1,3,4,5,7,\ldots)$ inductively by
$c_0=1$ and for $n\ge0$
\beq
\label{def1}
c_{n+1}=\case{c_n+1}{if $(c_n+1)/2\not\in\bc$,}{c_n+2}{otherwise.}
\eeq
Equivalently, $\bc$ is the lexicographically least sequence of
positive integers such that 
\beq
\label{def2}
\mbox{$m\in\bc$ if and only if $m/2\not\in\bc$.}
\eeq
It follows that $\bc$ contains all the positive odd
integers $m$ since in this case $m/2$ is not integral.  

The sequence $\bc$ is intimately connected with the {\it Thue-Morse
sequence\/} $\bt=(t_0,t_1,t_2,\ldots)=(0,1,1,0,1,0,0,1,\ldots)$ which is the
0-1 sequence defined inductively by $t_0=0$ and for $n\ge1$
$$
t_n=\case{t_{n/2}}{if $n$ even,}{1-t_{(n-1)/2}}{if $n$ odd.}
$$
A {\it block\/} of a sequence is a maximal subsequence of consecutive,
equal elements.  One can show~\cite{aab:srt} that $c_n-c_{n-1}$ is the
length of the $n$th block of $\bt$ (where we start with the 0th block
and set $c_{-1}=0$).

Given a sequence $\bs=(s_0,s_1,s_2,\ldots)$ and integers $k,l$ we
let
\beq
\label{lin}
k\bs+l=(ks_0+l,ks_1+l,ks_2+l,\ldots).
\eeq
To simplify our notation, we will also write  $k\Cong l\ (\Mod m)$ as
$k\Cong_m l$ with the usual convention that if the modulus is
omitted then $m=2$.
We can now characterize the parity of $M_n$.
\bth
\label{Mnthm}
The Motzkin number $M_n$ is even if and only if either  $n\in4\bc-2$
or $n\in4\bc-1$.
\eth
\pf
To prove this result we will need a combinatorial interpretation of
$M_n$.  A {\it 0-1-2 tree\/} is an ordered tree where each vertex has at
most two children (but a single child is not distinguished by being
either left or right).  It is known  that $M_n$ is
the number of 0-1-2 trees with $n$ edges.  
See the articles of Donaghey~\cite{don:rpt} and Donaghey and
Shapiro~\cite{ds:mn} for details.
The four 0-1-2 trees with
three edges are shown in Figure~\ref{012}

\thicklines
\setlength{\unitlength}{2pt}
\bfi
\bpi(40,60)(0,0)
\Gca \Gcc \Gce \Gcg
\Gcacc \Gccce \Gcecg
\epi
\hs{10pt}
\bpi(40,60)(0,0)
\Gac \Gec \Gce \Gcg
\Gacce \Gecce \Gcecg
\epi
\hs{30pt}
\bpi(40,60)(0,0)
\Gac \Gae \Gee \Gcg
\Gacae \Gaecg \Geecg
\epi
\hs{30pt}
\bpi(40,60)(0,0)
\Gec \Gae \Gee \Gcg
\Gecee \Gaecg \Geecg
\epi
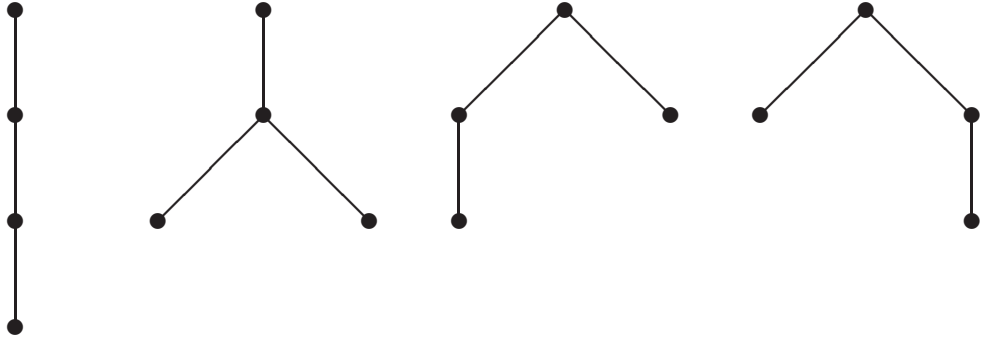
\capt{The 0-1-2 trees with three edges}\label{012}
\efi

Now let $S_n$ be the number of {\it symmetric 0-1-2 trees\/} which are
those with $n$ edges for which reflection in a vertical line
containing the root is an 
automorphism of the tree.  Only the first two trees in
Figure~\ref{012} are symmetric.  Clearly 
\beq
\label{MS}
M_n\Cong S_n
\eeq
for all $n\in\bbN$.  Furthermore, 
\beq
\label{S2n+1}
S_{2n+1}=S_{2n}
\eeq
since if a symmetric 0-1-2 tree has
$2n+1$ edges then the root must have a single child and the subtree
generated by that child must be a symmetric 0-1-2 tree with $2n$ edges.  So
to prove the theorem, it suffices to show that
$$
\mbox{$S_{2n}$ is even if and only if $2n\in4\bc-2$.}
$$
This can be restated that $S_{2n-2}$ is even iff $2n\in4\bc$ which is
equivalent to $n\in 2\bc$.  So, by~\ree{def2}, it suffices to prove
\beq
\label{S2n}
\mbox{$S_{2n-2}$ is even if and only if $n\not\in\bc$.}
\eeq

To prove~\ree{S2n}, we will need a recursion involving $S_{2n-2}$.  Let $T$ be
a symmetric 0-1-2 tree with $2n-2$ edges.  If the root of $T$ has
one child then the subtree generated by that child is a symmetric
0-1-2 tree with 
$2n-3$ edges.  If the root has two children then the subtree generated
by one
child can be any 0-1-2 tree with $n-2$ edges as long as the subtree
generated by the other is its reflection.  So using~\ree{MS} and~\ree{S2n+1}
\beq
\label{rr}
S_{2n-2}=S_{2n-3}+M_{n-2}\Cong S_{2n-4}+S_{n-2}.
\eeq

We now prove~\ree{S2n} by induction, where the case $n=1$ is trivial.
Suppose first that $n\not\in\bc$.  Then by~\ree{def1} we have $n-1\in\bc$ and
by induction it follows that  $S_{2n-4}=S_{2(n-1)-2}$ is odd.  Also, since
$n\not\in\bc$ we must have that $n$ is even.  Furthermore,
by~\ree{def2} we have $n/2\in\bc$.  By induction again,
$S_{n-2}=S_{2(n/2)-2}$ is odd.  So $S_{2n-4}+S_{n-2}$ is even and we
are done with this direction.

When $n\in\bc$, one can use similar reasoning to show that
$S_{2n-4}+S_{n-2}$ is odd.  One needs to consider the cases when $n$
is even and odd separately (and the latter case breaks into two
subcases depending on whether $n-1$ is in $\bc$ or not).  But there
are no really new ideas to the demonstration, so we omit the details.
\Qqed


We should note that Theorem~\ref{Mnthm} can also be derived from the
results in~\cite{kl:ipm}, although it is not explicitly stated there.
This theorem also permits us to determine the parity of various
related sequences which we will now proceed to do.

A {\it Motzkin path of length $n$\/}  is a lattice path in the lattice
$\bbN\times\bbN$ with steps $(1,1)$, $(1,-1)$, and $(1,0)$ starting at
$(0,0)$ and ending at $(n,0)$.  It is well known that $M_n$ is the
number of Motzkin paths of length $n$.  
(Note that we do not need any condition about staying above the
$x$-axis since we are working in $\bbN\times\bbN$.)
Define a {\it Motzkin prefix
of length $n$\/} to be a lattice path which forms the first $n$ steps
of a Motzkin path of length $m\ge n$.  Equivalently, a Motzkin prefix
is exactly like a Motzkin path except that the endpoint is not
specified.  Let $P_n$, $n\ge0$, be the number of Motzkin prefixes of length
$n$.  This is sequence A005773 in Sloane's Encyclopedia~\cite{slo:ole}.
The $P_n$ also count directed rooted animals with $n+1$ vertices as
proved by Gouyou-Beauchamps and Viennot~\cite{gbv:etd}.

\bco
\label{Pn}
The number $P_n$ is even if and only if $n\in2\bc-1$.
\eco
\pf
Let $s_n$ be the number of Motzkin paths of length $n$ which are
symmetric with respect to reflection in the vertical line $x=n/2$.
Clearly $M_n\Cong s_n$ for all $n\ge0$.  There is also a bijection
between Motzkin prefixes of length $n$ and symmetric Motzkin paths of
length $2n$ gotten by concatenating the prefix with its reflection in
the line $x=n$.  So $P_n=s_{2n}$.  Combining this with the previous
congruence and Theorem~\ref{Mnthm} completes the proof.
\Qqed

Next we consider the {\it Riordan numbers}~\cite[A005043]{slo:ole}
$\ga_n$ which count the number of ordered trees with $n$ edges where
every nonleaf has at least two children.  These are called {\it short
bushes\/} by Bernhart~\cite{ber:cmr}.  If we relax the degree
restriction so that the root can have any number of children then the
resulting trees are called {\it bushes}.  It is
known~\cite{don:rpt,ds:mn} that $M_n$ is the number of bushes with
$n+1$ edges.  It follows that
\beq
\label{Mg}
M_n=\ga_{n+1}+\ga_n
\eeq
since every bush with $n+1$ edges is either a short bush or has a root
with one child which generates a short bush with $n$ edges.

\bco
The number $\ga_n$ is even if and only if $n\in2\bc-1$.
\eco
\pf  Given the previous corollary, it suffices to show that $\ga_n$
and $P_n$ have the same parity.  So it suffices to show that the two
sequences satisfy the same recursion and boundary condition modulo 2.  
Now $\ga_0=1=P_0$ and we have just seen that 
$$
\ga_{n+1}\Cong \ga_n+M_n.
$$

So consider the prefixes $p$ counted by $P_{n+1}$.  If $p$ goes
through $(n,0)$ then there are two possible last steps for $p$ and so
such paths need not be considered modulo 2.  If $p$ goes through
$(n,m)$ where $m>0$ then those $p$ ending with a $(1,1)$ step can be
paired with those ending with a $(1,-1)$ step and ignored.  So we are
left with prefixes going through $(n,m)$ and $(n+1,m)$ where $m>0$.  Such
prefixes are equinumerous with those ending at $(n,m)$.  And since $m>0$,
this is precisely the set of Motzkin prefixes which are not Motzkin
paths.  So
$$
P_{n+1}\Cong P_n-M_n\Cong P_n+M_n
$$
as desired.
\Qqed

Finally, consider the sequence counting {\it restricted hexagonal
polyominos}~\cite[A002212]{slo:ole}.  The reader can find the precise
definition of these objects in the paper of Harary and
Read~\cite{hr:etp}.  We will use an equivalent definition in terms of
trees which can be obtained from the polyomino version by connecting
the centers of adjacent hexagons.  
A {\it ternary tree\/} is a rooted tree where every vertex has some subset
of three possible children: a left child, a middle child, or a right
child.  Just as with our definition of binary trees, this differs from
the all or none convention for ternary trees.
A {\it hex tree\/} is a ternary
tree  where no node can have two adjacent children.  (A middle child
would be adjacent to either a left or a right child but left and right
children are not adjacent.)  Let $H_n$, $n\ge0$, be the number of hex
trees having $n$ edges.

\bco
The number $H_n$ is even if and only if $n\in4\bc-2$ or $n\in4\bc-1$.
\eco
\pf
In view of Theorem~\ref{Mnthm}, it suffices to show that $H_n$ and
$M_n$ have the same parity.  Call a hex tree {\it symmetric\/} if the
reflection in a line containing the root leaves it invariant, and let
$h_n$ be the number of such trees with $n$ edges.  There is an obvious
bijection between symmetric hex trees and symmetric 0-1-2 trees.  So
$$
H_n\Cong h_n = S_n \Cong M_n
$$
as desired.
\Qqed

\section{Central binomial and trinomial coefficients}
\label{cbt}

Our main tool in this section will be the following famous congruence
of Lucas.  If the base $p$ expansion of $n$ is
$$
n=n_0+n_1 p +n_2 p^2 + \cdots
$$
then it will be convenient to denote the sequence of digits by
$$
(n)_p=(n_0,n_1,n_2,\ldots)=(n_i).
$$
\bth[Lucas~\cite{luc:cne}]
Let $p$ be a prime and let $(n)_p=(n_i)$ and
$(k)_p=(k_i)$.  Then
\beq
\label{luc}
{n\choose k} \Cong_p \prod_i {n_i\choose k_i}.\qqed
\eeq
\eth

The following corollary will be useful as well.  It is also a special
case of the theorem of Kummer cited in the discussion following the
statement of Theorem~\ref{Cnthm}.  But this
result will be sufficient for our purposes.
\bco
Let $p$ be prime.  If there is a carry when adding $k$ to $n-k$ in
base $p$ then 
$$
{n\choose k} \Cong_p 0.
$$
\eco
\pf
Using the notation of the previous theorem, if there is a carry out of
the $i$th place then we have $n_i<k_i$.  So ${n_i\choose k_i}=0$ and
thus the product side of~\ree{luc} is zero.
\Qqed

Most of our results in this section will have to do with congruences
modulo 3 so it will be useful to have the following notation.  Given
$i,j$ distinct integers in $\{0,1,2\}$ we let
$$
T(ij)=
\{n\in\bbN\ :\ \mbox{$(n)_3$ contains only digits equal to $i$ or $j$}\}.
$$

We begin with the central binomial coefficients. 
Recall that $\de_3(n)$ is the number of ones in the base three
expansion of $n$.  The next result
settles conjectures of Benoit Cloitre and Reinhard
Zumkeller~\cite[A074938--40]{slo:ole}.
\bth
\label{2ncn}
The central binomial coefficients satisfy
$$
{2n\choose n}\Cong_3 \case{(-1)^{\de_3(n)}}{if $n\in T(01)$,}{0}{otherwise.}
$$
\eth
\pf  
If $n$ has a 2 in its ternary expansion then there is a carry when
adding $(n)_3$ to itself.  So the second half of the theorem
follows from the previous corollary.  On the other hand, if 
$n\in T(01)$ then $2n\in T(02)$ and $(2n)_3$ has twos exactly where
$(n)_3$ has ones.  So by Lucas' Theorem
$$
{2n\choose n}\Cong_3 {2\choose 1}^{\de_3(n)} \Cong_3 (-1)^{\de_3(n)} 
$$ 
giving the first half.
\Qqed

It is easy to generalize the previous theorem to arbitrary prime
modulus.  
To state the result, we need to define
\beq
\label{depj}
\de_{p,j}(n)=\mbox{number of elements of $(n)_p$ equal to $j$}
\eeq
where $0\le j<p$.  
Since the proof of the general case is the same as the one just given, we omit it.
\bth
\label{2ncnp}
Let $p$ be prime and let $S$ be the set of all $n\in\bbN$ such that all
elements of $(n)_p$ are less than  or equal to $p/2$.  Then
$$
{2n\choose n}\Cong_p
\case{\dil\prod_j {2j\choose j}^{\de_{p,j}(n)}}{if $n\in S$,}
{0}{otherwise.\rule{0pt}{20pt}\qqed}
$$
\eth

It turns out that there is a connection between the central binomial
coefficients modulo 3 and the Thue-Morse sequence $\bt$.  This may
seem surprising because $\bt$ is essentially a modulo 2 object.
However, Theorem~\ref{2ncn} will allow us to reduce questions about
${2n\choose n}$ mod 3 to questions about bit strings.  We will need
another one of the many definitions of $\bt$ for the proof, namely
\beq
\label{tn}
t_n=\rho(\de(n))
\eeq
where $\rho(k)$ is the remainder of $k$ on division by 2.
We will also need the notation that $\ba\Cong_m \bb$ as sequences if
and only if $a_n\Cong_m b_n$ for all $n\in\bbN$.
The next result is again a conjecture of
Cloitre~\cite[A074939]{slo:ole}.
\bth
\label{2ncn1}
We have
$$
\left(n\ :\ {2n\choose n}\Cong_3 1 \right) \Cong_3 \bt.
$$
\eth
\pf
Let us call the sequence on the left of the previous congruence
$\ba$.  Then from Theorem~\ref{2ncn} we have that $n\in\ba$ exactly
when $n\in T(01)$ and $(n)_3$ has an even number of ones.  From this
it follows by an easy induction that $n=a_i$ if and only if
$(n)_3=(n_0,n_1,n_2,\ldots)$ where $(i)_2=(n_1,n_2,\ldots)$ and $n_0$
is zero or one depending on whether $\de(i)$ is even or odd,
respectively.  So by~\ree{tn} we have
$$
a_i=n\Cong_3 n_0=\rho(\de(i))=t_i
$$
for all $i\ge0$.
\Qqed

There is an analogous conjecture of Cloitre for those central binomial
coefficients with residue $-1$ modulo 3~\cite[A074938]{slo:ole}.  Since
the proof is much the same as the previous one, we omit it.
\bth
We have
$$
\left(n\ :\ {2n\choose n}\Cong_3 -1 \right) \Cong_3 1-\bt.\Qqed
$$
\eth

We next consider the central trinomial
coefficients~\cite[A002426]{slo:ole}.  Let $T_n$ be the largest
coefficient in the expansion of $(1+x+x^2)^n$.  It is
easy~\cite{bps:mf} to express $T_n$ in terms of trinomial coefficients
\beq
\label{Tn}
T_n=\sum_{k\ge0} {n\choose k,\ k,\ n-2k}
\eeq
where we use the convention that if any multinomial coefficient has a
negative number on the bottom then the coefficient is zero.  Lucas'
Theorem and its corollary generalize in the expected way to
multinomial coefficients.  So now we can find the residue of $T_n$ modulo 3.
\bth
\label{Tnthm}
The central trinomial coefficients satisfy
$$
T_n\Cong_3
\case{1}{if $n\in T(01)$,}{0}{otherwise.}
$$
\eth
\pf
Modulo 3 we can restrict the sum in~\ree{Tn} to 
those $k$ such that there is no
carry in doing the triple addition $k+k+(n-2k)$ in base 3.  So, in
particular, we can restrict to $k\in T(01)$ since
if $(k)_3=(k_i)$ contains a 2 then we will have such a
carry.  Furthermore, if
$k_i=1$ for some $i$ then $k+k$ has a two in the $i$th place, and to have
no carry this forces $n_i=2$.

Now let $(n-2k)_3=(l_i)$ and let $S$ be the set of indices
$i$ such that $n_i=2$.  So we have shown that $\De_3(k)\sbe S$.
Furthermore, for every $i\not\in\De_3(k)$ we must have $l_i=n_i$ since
$k\in T(01)$.  So the nonzero terms in the sum correspond to subsets
$R\sbe S$ and each such subset contributes
$$
{2\choose 1,\ 1,\ 0}^{|R|}=2^{|R|}.
$$
Hence, by the binomial theorem, the total contribution is
$$
\sum_{R\sbe S} 2^{|R|}=3^{|S|}\Cong_3
\case{1}{if $S=\emp$,}{0}{if $S\neq\emp$.}
$$
But $S=\emp$ precisely when $n\in T(01)$, so we are done.
\Qqed

Since the $T_n$ are related to a number of the other sequences which
we have been studying, we can use the previous result to determine
their behavior modulo 3.  We will apply linear operations to sets the
same way we do to sequences~\ree{lin}.
\bco
The Motzkin numbers satisfy
$$
M_n\Cong_3
\left\{\barr{cl}
-1&\mbox{if $n\in 3T(01)-1$,}\\
1&\mbox{if $n\in 3T(01)$ or $n\in 3T(01)-2$,}\\
0&\mbox{otherwise.}
\earr\right.
$$
\eco
\pf
Barcucci, Pinzani, and Sprugnoli~\cite{bps:mf} have shown that
\beq
\label{MT}
2M_n=3T_n+2T_{n+1}-T_{n+2}.\qquad 
\eeq
Reducing this equation modulo 3 and applying the previous theorem
finishes the proof.
\Qqed

\bco
The Motzkin prefix numbers satisfy
$$
P_n\Cong_3
\left\{\barr{cl}
1&\mbox{if $n\in 3T(01)$,}\\
-1&\mbox{if $n\in 3T(01)+1$ or $3T(01)-1$,}\\
0&\mbox{otherwise.}
\earr\right.
$$
\eco
\pf
If $p$ is a Motzkin prefix of length $n$ going through $(n-1,m)$ for
some $m>0$ then there are three ways to end the prefix and so they
cancel out modulo 3.  If $p$ goes through $(n-1,0)$ then the first
$n-1$ steps of $p$ form a Motzkin path and there are two possible last
steps.  So $P_n\Cong_3 2M_{n-1}$.  Now apply the previous corollary to
finish.
\Qqed

\bco
The Riordan numbers satisfy
$$
\ga_n\Cong_3
\case{1}{if $n\in T(01)-1$,}{0}{otherwise.}
$$
\eco
\pf
Using the recursions~\ree{Mg} and~\ree{MT} it is easy to prove
inductively that $\ga_n\Cong_3 T_{n+1}$.  Theorem~\ref{Tnthm} now
completes the proof.
\Qqed

\section{Miscellaneous results and conjectures}
\label{mrc}

We end with various results and conjectures related to what we have
done in the previous sections.

\subsection{Catalan numbers to other moduli}

Theorem~\ref{Cnthm} implies that the $k$th block of
zeros in the sequence of Catalan numbers modulo 2 has length $2^k-1$
(where we start numbering with the first block).  Alter and
Kubota~\cite{ak:ppp} have generalized this result to arbitrary primes
and prime powers.  One of their main theorems is as follows.
\bth(Alter and Kubota)
\label{ak}
Let $p\ge3$ be a prime and let $q=(p+1)/2$. The length of the $k$th
block of zeros of the Catalan sequence modulo $p$ is
$$
\left(p^{\om_q(k)+\de_{3,p}+1}-3\right)/2
$$
where $\de_{3,p}$ is the Kronecker delta.\Qqed
\eth

We can improve on this theorem in several regards.  First of all, when
$p=3$ we can use our results to give a complete characterization of
the residue of $C_n$ and not just say when it is zero.  Suppose
$(n)_3=(n_i)$.  Then we let
$$
T^*(01)=\{n\ :\ \mbox{$n_i=0$ or 1 for all $i\ge1$}\}
$$
and
$$
\de_3^*(n)=\mbox{number of $n_i=1$ for $i\ge1$.} 
$$
\bth
\label{Cn3}
The Catalan numbers satisfy
$$
C_n\Cong_3
\case{(-1)^{\de_3^*(n+1)}}{if $n\in T^*(01)-1$,}{0}{otherwise.}
$$
\eth
\pf
The result is easy to verify for $n\le1$ so we assume $n\ge2$.
Directly from our definition of $C_n$ we have
$$
C_n=\frac{4n-2}{n+1}C_{n-1}
$$
If $n\Cong_3 0$ or $1$ then $n+1$ is invertible modulo 3 and in fact
$(4n-2)/(n+1)\Cong_3 1$.  So for $k\ge 1$ we have 
$C_{3k-1}\Cong_3 C_{3k}\Cong_3 C_{3k+1}$.  Thus it suffices to prove
the theorem for $n\Cong_3 0$.  Notice that in this case 
$C_n\Cong_3 {2n\choose n}$.  Furthermore
$n+1\in T^*(01)$ if and only if $n\in T(01)$.  And lastly
$\de_3^*(n+1)=\de_3(n)$.  Applying Theorem~\ref{Cnthm} finishes the proof.
\Qqed

We should verify that we can derive the $p=3$ block lengths in
Theorem~\ref{ak} from Theorem~\ref{Cn3}.
First from the latter result it follows that the $k$th block must
start at an integer $3a-1$ and end at $3b-1$ for $a,b\in\bbN$.  To
simplify notation, let $\om=\om_2(k)$.
Now $(a)_3$ must contain a 
2 and $(a-1)_3$ does not.  It follows that
$(a)_3=(a_0,a_1,a_2,\ldots)$ where 
$a_0=2$ and $(a_1,a_2,\ldots)=(k-1)_2$.
Furthermore, since $b+1$ is the smallest integer larger than $a$ whose
expansion contains no twos, the first $\om+1$ elements of $(b)_3$ must
all equal 2 and the rest must agree with the corresponding entries of
$(a)_3$.  By the same token, we must have
$a_1=a_2=\ldots=a_\om=1$. Now one calculates the number of
integers in the $k$th block by considering the first $\om+1$
digits of $a$ and $b$ to get a count of
$$
3(a-b+1)=3[(3^{\om+1}-1)-(3^\om+3^{\om-1}+\cdots+3+2)+1]=(3^{\om+2}-3)/2
$$
as desired.  Note that not only have we been able to determine the
length and starting and ending points of the block (which was also done
by Alter and Kubota) but our demonstration is combinatorial as opposed
to the original proof of Theorem~\ref{ak} which is arithmetic.  We had
to use Lucas' Theorem to get to this result, but that theorem also has
a combinatorial demonstration using group actions~\cite{sag:cva}.  

When $p\ge5$, the residues of $C_n$ become more complicated, but one
could use the same techniques in principle to compute them.  In
particular, if one is only interested in divisibility then one can
derive Theorem~\ref{ak} from Theorem~\ref{2ncnp} as we did for the
$p=3$ case above. 

It is also interesting another setting where a congruence involving
the Catalan numbers has arisen.  Albert, Atkinson, and
Klazar~\cite{aak:esp} have studied {\it simple permutations\/} which
are those permutations of $\{1,2,\ldots,n\}$ mapping no nontrivial
subinterval of this set onto an interval.  Then the number of such simple
permutations is $2(-1)^{n+1}-\Com_n$ where $\Com_n$ is the
coefficient of $x^n$ in the compositional inverse of the formal power
series $\sum_{n\ge1} n!x^n$~\cite[A059372]{slo:ole}.  One of the
results in~\cite{aak:esp} is that
$$
\Com_n\Cong_3 C_{n-1}.
$$
Their proof of this result uses generating functions, so it would be
interesting to find a combinatorial one.  Also, one would like to know
the behavior of $\Com_n$ modulo other odd primes.  (Albert et.\ al.\
have results for powers of two.)

The careful reader will note that we have not yet derived the
residues of the hex tree numbers $H_n$ modulo three.  It is
time to fill that lacuna.
\bth
The hex tree numbers satisfy
$$
H_n\Cong_3
\case{(-1)^{\de_3^*(m+1)}}{if $n=2m$ where $m\in T^*(01)-1$,}{0}{otherwise.}
$$
\eth
\pf
Suppose $T$ is a hex tree which has a vertex with a single child.  Finding
the first such vertex, say in depth-first order, one can associate
with $T$ the two other hex trees which differ from $T$ only by moving
the child into the two other possible positions.  So modulo 3, $H_n$
is congruent to the number of hex tree with $n$ edges where every
vertex has 0 or 2 children.  So to be nonzero modulo 3, we must have
$n=2m$ and the resulting trees are in bijection with binary trees on
$m$ vertices (merely remove the $m$ leaves of the hex tree).  
Thus $H_n\Cong_3 C_m$ and
we are now done by Theorem~\ref{Cn3}.
\Qqed

\subsection{Motzkin numbers to other moduli}

For the Motzkin numbers, nothing has been proved for moduli other
than 2 and 3.  However, there are some conjectures.  
To put them in the context of Theorem~\ref{Mnthm}, note that the
Thue-Morse block sequence $\bc$ can also be described~\cite{aab:srt}
as the increasing sequence of all numbers of the form
$$
\mbox{$(2i+1)4^j$ where $i,j\in\bbN$.}
$$
The following conjecture is due in part to Tewodros
Amdeberhan~\cite{amd:pc}. 
\bcon[Amdeberhan]
We have $M_n\Cong_4 0$ if and only if
$$
\mbox{$n=(4i+1)4^{j+1}-1$ or $n=(4i+3)4^{j+1}-2$ where $i,j\in\bbN$.}
$$
Furthermore we never have $M_n\Cong_8 0$.
\econ
Amdeberhan also has a conjecture about some of the values of $n$ for which
$M_n$ is zero modulo 5, although it is complicated.

\subsection{Gould's sequence}

{\it Gould's sequence}~\cite[A001316]{slo:ole} consists of the numbers
$G_n$ which count the number of odd entries in the $n$th row of
Pascal's triangle.  More generally, we can calculate $G_n(p)$ which is
the number of entries in the $n$th row of Pascal's triangle which are
not zero modulo the prime $p$.  Recall the definition of
$\de_{p,j}(n)$ in~\ree{depj}. 
\bth
Let $p$ be prime.  Then
$$
G_n(p)=\prod_{1\le j<p} (j+1)^{\de_{p,j}(n)}.
$$
Furthermore, every entry of the $n$th row of Pascal's triangle is
nonzero modulo $p$ if and only if
$$
n=qp^k-1
$$
where $1\le q<p$ and $k\in\bbN$.  In particular
$$
G_n=2^{\de(n)}
$$
and every entry of the $n$th row of Pascal's triangle is odd if and
only if $n=2^k-1$ where $k\in\bbN$.
\eth
\pf
Suppose ${n\choose k}\not\Cong_p 0$ where $(n)_p=(n_i)$ and $(k)_p=(k_i)$.  
If $n_i=j$ then we will not have a carry in the $i$th place if and
only if $0\le k_i\le j$.  So there are $j+1$ choices for $k_i$ and
taking the product of the number of choices for each $i$ gives the
first statement of the theorem.

Now suppose that every entry of the $n$th row is nonzero modulo $p$.
Since there are no carries for all $k$, all the elements of $(n)_p$
must equal $p-1$ except for possibly the last (leading) one $n_l$.
Since there can never be a carry out of $n$'s last place, we have the
desired characterization of those $n$ under consideration.
\Qqed

\subsection{Sums of central binomial coefficients}
The partial sums of central binomial
coefficients~\cite[A006134]{slo:ole} also have nice
congruence properties.  The proof of the next result is easily
obtained by using Theorem~\ref{2ncn} and induction on $n$, so we omit
it.  In conjunction with Theorem~\ref{2ncn1}, it settles a conjecture
of Cloitre~\cite[A083096]{slo:ole}.
\bth
We have
$$
\sum_{k\ge0} {2k\choose k}\Cong_3
\case{(-1)^{\de_3(n)}}{if $n\in 3T(01)$,}{0}{otherwise.\Qqed}
$$
\eth

\subsection{Ap\'ery numbers and central Delannoy numbers}

We can generalize our results about the central trinomial numbers as
follows.  Given positive integers $r,s$ we define a sequence with the
following entries
$$
a_n(r,s)=\sum_{k\ge0} {n\choose k}^r {n+k\choose k}^s.
$$
Note that since $r,s$ are positive, each term in this sum will have a
factor of
$$
{n\choose k}{n+k\choose k}={n+k\choose k,\ k,\ n-k}.
$$
Using this fact we can prove the following result.  Since the
demonstration is similar to that of Theorem~\ref{Tnthm}, it is
omitted.  Again, this settles a conjecture of Cloitre~\cite{clo:pc}.
\bth
Let $r,s$ be positive integers.  Then
$$
a_n(r,s)\Cong_3
\left\{\barr{cl}
(-1)^{\de_3(n)}&\mbox{if $s$ is even,}\\
1&\mbox{if $s$ is odd and  $n\in T(02)$,}\\
0&\mbox{otherwise.\Qqed}
\earr\right.
$$
\eth

The {\it central Delannoy numbers}~\cite[A001850]{slo:ole} are
$D_n=a_n(1,1)$.  Also, the 
{\it Ap\'ery numbers}~\cite[A005258]{slo:ole} are 
$A_n=a_n(2,1)$.  So we immediately have the following corollary.
\bth
The central Delannoy numbers and Ap\'ery numbers satisfy
$$
D_n\Cong_3 A_n\Cong_3
\case{1}{if $n\in T(02)$,}{0}{otherwise.\Qqed}
$$
\eth 

\subsection{Central Eulerian numbers}

The {\it Eulerian numbers}~\cite[A008292]{slo:ole} are denoted
$A(n,k)$ and count the number of permutations in the symmetric group
$\fS_n$ which have $k-1$ descents.  They can be written as
$$
A(n,k)=\sum_{i=0}^k (-1)^i (k-i)^n {n+1\choose i}.
$$
Since the odd numbered rows have an odd number of elements, we define
the {\it central Eulerian numbers\/} to be
$$
E_n=A(2n-1,n)=\sum_{i=0}^n (-1)^i (n-i)^{2n-1}{2n\choose i}.
$$
We have the following congruence for these numbers.
\bth
The central Eulerian numbers satisfy
$$
E_n\Cong_3
\case{1}{if $n\in T(01)+1$,}{0}{otherwise.}
$$
\eth
\pf
Note that $k^{2n-1}=k$ for $k=0,\pm 1$.  So we have
$$
E_n\Cong_3 \sum_{i=0}^n (-1)^i (n-i){2n\choose i}.
$$
Applying the binomial recursion to this sum twice yields, after
massive cancellation,
$$
E_n\Cong_3 (-1)^{n-1}{2n-2\choose n-1}.
$$
Now Theorem~\ref{2ncn} will finish the proof provided $n+\de_3(n)$ is
always even.  But this is easy to show by induction on $n$, so we are
done.
\Qqed

Rows in the Eulerian triangle are symmetric, so even numbered rows
have two equal elements in the middle.  We will call these elements
{\it bicentral}.  Cloitre conjectured the residues of these elements
modulo 3.  Since the proof of this result is similar to the one just
given, we will omit it.
\bth
The bicentral Eulerian numbers satisfy
$$
A(2n,n)\Cong_3
\left\{\barr{cl}
1&\mbox{if $n\in 3T(01)+1$,}\\
-1&\mbox{if $n\in 3T(01)$ or $3T(01)+2$,}\\
0&\mbox{otherwise.\Qqed}
\earr\right.
$$
\eth

\subsection{Noncrossing connected graphs}

Noncrossing set partitions are an important object of study in
combinatorics.  An excellent survey of the area can be found in the
article of Simion~\cite{sim:np}.  Noncrossing graphs are a
generalization of noncrossing partitions which have been studied by
Flajolet and Noy~\cite{fn:acn}.  Consider vertices labeled
$1,\ldots,n$ and arranged clockwise in this order around a circle.
A graph on this vertex set is  {\it noncrossing\/} if, when the edges
are drawn with straight line segments between the vertices, none of the
edges cross.  Let $N_n$ be the number on noncrossing connected graphs
on $n$ vertices~\cite[A007287]{slo:ole}.  Then it can be shown that
$$
N_n=\frac{1}{n-1}\sum_{k\ge0} {3n-3\choose n+k+1}{k\choose n-2}.
$$
We have the following conjecture about the residue of $N_n$ modulo 3.
\bcon
The number of noncrossing connected graphs satisfies
$$
N_n\Cong_3
\left\{\barr{cl}
1&\mbox{if $n=3^i$ or $n=2\cdot3^i$ for some $i\in\bbN$,}\\
-1&\mbox{if $n=3^i+3^j$ for two distinct $i,j\in\bbN$,}\\
0&\mbox{otherwise.}
\earr\right.
$$
\econ

In the first two cases, it is not hard to show that the congruence
holds using Lucas' Theorem  because of the very specific form of 
$(n)_3$.  However, we have been unable to prove that for all remaining
$n$ one always has $N_n$ divisible by $3$.  It would be even more
interesting to give a combinatorial proof of this result based on
symmetries of the graphs involved.

{\bf Acknowledgements}  We would like to thank Neil Sloane for
maintaining the On-Line Encyclopedia of Integer Sequences.  It has
been invaluable to us in helping to form conjectures about various
sequences.

\begin{\bib}{99}

\bibitem{aak:esp} M. Albert, M. Atkinson, and M. Klazar,
The enumeration of simple permutations,
{\it J. Integer Sequences\/} {\bf 6} (2003), Article 03.4.4, 18 pp.

\bibitem{aab:srt} J.-P. Allouche, A. Arnold, J. Berstel, S. Brlek, 
W. Jockusch, S. Plouffe, and B. E. Sagan, A relative of the
Thue-Morse sequence, {\it Discrete Math.\/} {\bf 139} (1995), 455--461.

\bibitem{as:ups} J.-P. Allouche and J. Shallit, The ubiquitous
Prouhet-Thue-Morse sequence, in ``Sequences and their applications,
Proceedings of SETA'98,'' C. Ding, T. Helleseth, and H. Niederreiter
eds., Springer-Verlag, 1999, 1--16.

\bibitem{ak:ppp} R. Alter and K. Kubota, Prime and prime power
divisibility of Catalan numbers,
{\it J. Combin.\ Theory Ser.\ A} {\bf 15} (1973), 243--256.

\bibitem{amd:pc} T. Amdeberhan, personal communication.

\bibitem{bps:mf} E. Barcucci, R. Pinzani, R. Sprugnoli, The Motzkin
family, {\it Pure Math.\ Appl.\ Ser.\ A\/} {\bf 2} No. 3--4 (1991), 249--279.

\bibitem{ber:cmr} F. Bernhart, Catalan, Motzkin, and Riordan numbers,
{\it Discrete Math.\/} {\bf 204} (1999), 73--112.

\bibitem{clo:pc} B. Cloitre, personal communication.

\bibitem{deu:idp}  E. Deutsch, An involution on Dyck paths and its
consequences, {\it Discrete Math.\/} {\bf 204} (1999), 163--166.

\bibitem{dic:htn}  L. E. Dickson, ``History of the Theory of Numbers,
Vol.\ 1,''  Chelsea, New York, NY, 1952.

\bibitem{don:rpt} R. Donaghey, Restricted plane tree representations
of four Motzkin-Catalan equations,
{\it J. Combin.\ Theory, Ser.\ B\/} {\bf 22} (1977), 114--121.

\bibitem{ds:mn} R. Donaghey and L. W. Shapiro, Motzkin numbers,
{\it J. Combin.\ Theory, Ser.\ A\/} {\bf 23} (1977), 291--301.

\bibitem{ege:pcn} \"O. E\u{g}ecio\u{g}lu, The parity of the Catalan
numbers via lattice paths, {\it Fibonacci Quart.\/}  {\bf 21} (1983)
65--66. 

\bibitem{fn:acn} P. Flajolet and M. Noy, Analytic combinatorics of
non-crossing configurations, {\it Discrete Math.\/} {\bf 204} (1999),
203--229.

\bibitem{gbv:etd} D. Gouyou-Beauchamps and G. Viennot, Equivalence of
the two-dimensional directed animal problem to a one-dimensional path
problem, {\it Adv.\ Appl.\ Math.\/} {\bf 9} (1988), 334--357.

\bibitem{hr:etp} F. Harary and R. C. Read, The enumeration of
tree-like polyhexes, {\it Proc.\ Edinburgh Math.\ Soc.\/} {\bf 17}
(1970), 1--13.

\bibitem{kl:ipm} M. Klazar and F. Luca, On integrality and periodicity
of the Motzkin numbers, preprint.

\bibitem{kum:ear} E. E. Kummer \"Uber die Erg\"anzungss\"atze zu den
allgemeinen Reciproci\-t\"atsgesetzen, {\it J. Reine Angew.\ Math.\/}
{\bf 44} (1852) 93--146.

\bibitem{luc:pfm} F. Luca, Prime factors of Motzkin numbers, preprint.

\bibitem{luc:cne} E. Lucas, Sur les congruences des nombres eul\'eriens et
des coefficients diff\'erentiels des fonctions trigonom\'etriques
suivant un module premier, {\it Bull.\ Soc.\ Math.\ France\/} {\bf 6}
(1877--1878), 49--54.


\bibitem{sag:cva} B. E. Sagan, Congruences via Abelian groups,
{\it J. Number Theory} {\bf 20} (1985), 210--237.

\bibitem{sch:vcp} E. Schr\"oder, Vier combinatorische Probleme, 
{\it Z. f\"ur Math.\ Phys.\/} {\bf 15} (1870), 361--376.

\bibitem{sim:np} R. Simion, Noncrossing partitions,
{\it Discrete Math.\ } {\bf 217} (2000), 367--409.

\bibitem{su:sln} R. Simion and D. Ullman,  On the structure of the
lattice of noncrossing partitions, {\it Discrete Math.\/} {\bf 98}
(1991), 193--206.

\bibitem{slo:ole} N. J. A. Sloane, ``The On-Line Encyclopedia of
Integer Sequences,'' available at
{\bf http://www.research.att.com/\~{\rule{1pt}{0pt}}njas/sequences/}.

\bibitem{sta:ec2} R. P. Stanley, ``Enumerative Combinatorics,
Volume 2,''  Cambridge University Press, Cambridge, 1999.

\end{\bib}

\end{document}